\documentclass[11pt,a4paper]{amsart}

\usepackage{amssymb}
\usepackage{rotating}

\newcommand{\1}{{1\!\!1}}
\usepackage{array}
\newcolumntype{C}{>{$}c<{$}}
\newcolumntype{L}{>{$}l<{$}}
\newcolumntype{R}{>{$}r<{$}}

\def\to{\mathchoice
{\longrightarrow}
{\rightarrow}
{\rightarrow}
{\rightarrow}}

\def\mapsto{\DOTSB\mapstochar\to}

\newtheorem{theorem}[subsection]{Theorem}
\newtheorem*{theorem*}{Theorem}
\newtheorem{proposition}[subsection]{Proposition}
\newtheorem{lemma}[subsection]{Lemma}
\newtheorem{corollary}[subsection]{Corollary}

\theoremstyle{definition}

\newtheorem{definition}[subsection]{Definition}
\newtheorem{example}[subsection]{Example}
\newtheorem*{remark}{Remark}

\makeatletter\let\c@equation=\c@subsection\makeatother

\newcommand{\Z}{\mathbb{Z}}
\newcommand{\C}{\mathbb{C}}

\newcommand{\Q}{\mathbb{Q}}

\newcommand{\om}{\omega}

\renewcommand{\o}{\otimes}
\newcommand{\p}{\partial}

\newcommand{\half}{\tfrac12}

\newcommand{\bull}{\bullet}
\newcommand{\ohat}{\Hat{\otimes}}

\DeclareMathOperator{\Aut}{Aut}
\DeclareMathOperator{\Ind}{Ind}
\DeclareMathOperator{\Res}{Res}
\DeclareMathOperator{\res}{res}
\DeclareMathOperator{\Tr}{Tr}
\DeclareMathOperator{\rk}{rk}

\renewcommand{\]}{{]\!]}}
\renewcommand{\[}{{[\![}}
\renewcommand{\(}{(\!(}
\renewcommand{\)}{)\!)}

\DeclareMathOperator{\SL}{SL}

\DeclareMathOperator*{\colim}{colim}
\DeclareMathOperator{\Ob}{Ob}
\DeclareMathOperator{\ch}{ch}
\DeclareMathOperator{\Ch}{Ch}
\DeclareMathOperator{\Exp}{Exp}
\DeclareMathOperator{\Log}{Log}

\newcommand{\Cat}{{\mathcal C}}
\newcommand{\CM}{\mathcal{M}}
\newcommand{\Mbar}{\overline{\mathcal{M}}}
\renewcommand{\SS}{\mathbb{S}}
\newcommand{\SSS}{\mathsf{S}}
\DeclareMathOperator{\Serre}{\mathsf{e}}
\newcommand{\LL}{\mathsf{L}}

\DeclareMathOperator{\VERT}{Vert}
\newcommand{\V}{\mathcal{V}}
\newcommand{\w}{\mathcal{W}}
\newcommand{\CL}{\mathcal{L}}
\newcommand{\HH}{\mathsf{H}}
\newcommand{\CO}{\mathcal{O}}
\newcommand{\M}{\mathsf{M}}
\newcommand{\MM}{\mathbb{M}}
\newcommand{\KHM}{{\mathsf{KHM}}}

\renewcommand{\AA}{\mathbf{a}}
\newcommand{\BB}{\mathbf{b}}

\begin{document}

\title{The semi-classical approximation for modular operads}

\author{E. Getzler}

\address{Max-Planck-Institut f\"ur Mathematik, Gottfried-Claren-Str.\ 26,
D-53225 Bonn, Germany}

\curraddr{Department of Mathematics, Northwestern University, Evanston, IL
60208-2730, USA}

\email{getzler@math.nwu.edu}

\maketitle

The semi-classical approximation is an explicit formula of mathematical
physics for the sum of Feynman diagrams with a single circuit. In this
paper, we study the same problem in the setting of modular operads
\cite{modular}; instead of being a number, the interaction at a vertex of
valence $n$ will be an $\SS_n$-module.

The motivation for developing this theory was the desire to calculate the
$\SS_n$-equivariant Hodge polynomials of the Deligne-Mumford-Knudsen moduli
spaces $\Mbar_{1,n}$ of stable curves of genus $1$ with $n$ marked smooth
points. In performing these calculations, we use the formulas for the
$\SS_n$-equivariant Serre polynomials of $\Mbar_{0,n}$ and $\Mbar_{1,n}$
derived in \cite{gravity} and \cite{II} respectively.

A particular consequence of our calculations will be needed in
\cite{genus1} to find a relation among the codimension two cycles in
$\Mbar_{1,4}$.
\begin{theorem*}
The $\SS_4$-module $H^4(\Mbar_{1,4},\Q)$ is isomorphic to
$$
\bigl( V_{(4)}\o\Q^7 \bigr) \oplus \bigl( V_{(3,1)}\o\Q^4 \bigr) \oplus
\bigl( V_{(2,2)}\o\Q^2 \bigr) .
$$
\end{theorem*}

\subsection*{Acknowledgments}

I wish to thank the Department of Mathematics at the Universit\'e de
Paris-VII the Max-Planck-Institut f\"ur Mathematik in Bonn for their
hospitality during the inception and completion, respectively, of this
paper. I am grateful to D. Zagier for showing me the asymptotic expansion
of Corollary \ref{don}.

This research was partially supported by a research grant of the NSF and a
fellowship of the A.P. Sloan Foundation.

\section{Wick's theorem and the semi-classical approximation}

Let $\Gamma_{g,n}$ be the small category whose objects are isomorphism
classes of stable graphs $G$ of genus $g(G)=g$ with $n$ totally ordered
legs \cite{modular}, and whose morphisms are the automorphisms: if
$G\in\Gamma_{g,n}$, its automorphism group $\Aut(G)$ is the subset of the
permutations of the flags which preserve all the data defining the stable
graph, including the total ordering of the legs. Because of the stability
condition, $\Gamma_{g,n}$ is a finite category.

Define polynomials $\{\M v_{g,n}\mid 2(g-1)+n>0\}$ of a set of variables
$\{v_{g,n}\mid 2(g-1)+n>0\}$ by the following formula:
\begin{equation} \label{M}
\M v_{g,n} = \sum_{G\in\Ob\Gamma_{g,n}} \frac{1}{|\Aut(G)|}
\prod_{v\in\VERT(G)} v_{g(v),n(v)} .
\end{equation}
Introduce the sequences of generating functions
$$
a_g(x) = \sum_{2(g-1)+n>0} v_{g,n} \frac{x^n}{n!} , \quad\text{and}\quad
b_g(x) = \sum_{2(g-1)+n>0} \M v_{g,n} \frac{x^n}{n!} .
$$
Wick's theorem gives an integral formula for the generating functions
$\{b_g\}$ in terms of $\{a_g\}$:
$$
\sum_{g=0}^\infty b_g \hbar^{g-1} = \log \int_{-\infty}^\infty \exp \biggl(
\sum_{g=0}^\infty a_g \hbar^{g-1} - \frac{(x-\xi)^2}{2\hbar} \biggr) \,
\frac{dx}{\sqrt{2\pi\hbar}} .
$$
As written, this is purely formal, since it involves the integration of a
power series in $x$. It may be made rigourous by observing that the
integral transform
$$
f \mapsto \int_{-\infty}^\infty f(\hbar,x) e^{-(x-\xi)^2/2\hbar} \,
\frac{dx}{\sqrt{2\pi\hbar}}
$$
induces a continuous linear map on the space of Laurent series
$\Q\(\hbar\)\[x\]$ topologized by the powers of the ideal $(\hbar,x)$.

The semi-classical expansion is a pair of formulas for $b_0$ and
$b_1$ in terms of $a_0$ and $a_1$, which we now recall.
\begin{definition}
Let $R$ be a ring of characteristic zero. The Legendre transform $\CL$ is
the involution of the set $x^2/2+x^3R\[x\]$ characterized by the formula
$$
(\CL f)\circ f' + f = p_1 f' .
$$
\end{definition}
\begin{theorem} \label{legendre}
The series $x^2/2+b_0$ is the Legendre transform of $x^2/2-a_0$.
\end{theorem}

The first few coefficients of $b_0$ may be calculated, either from the
definition of $\M v_{0,n}$ or from Theorem \ref{legendre}:
$$\begin{tabular}{|C|L|} \hline
n & \M v_{0,n} \\[2pt] \hline
3 & v_{0,3} \\[5pt]
4 & v_{0,4} + 3v_{0,3}^2 \\[5pt]
5 & v_{0,5} + 10v_{0,4}v_{0,3} + 15v_{0,3}^3 \\[5pt]
6 & v_{0,6} + 15v_{0,5}v_{0,3} + 10v_{0,4}^2 + 105 v_{0,4}v_{0,3}^2 +
105 v_{0,3}^4 \\[3pt] \hline
\end{tabular}$$

We now come to the formula for $b_1$, known as the semi-classical
approximation.
\begin{theorem} \label{semi}
The series $b_1$ and $a_1$ are related by the formula
$$
b_1 = \bigl( a_1 - \tfrac{1}{2} \log (1-a_0'') \bigr) \circ ( x + b_0' ) .
$$
\end{theorem}

By the definition of the Legendre transform, we see that $(\CL f)'\circ
f'=x$. It follows that Theorem \ref{semi} is equivalent to the formula
$$
b_1 \circ ( x - a_0' ) =  a_1 - \tfrac{1}{2} \log(1-a_0'') .
$$
This formula expresses the fact that the stable graphs contributing to
$b_1$ are obtained by attaching a forest whose vertices have genus $0$
to two types of graphs:
\begin{enumerate}
\item those with a single vertex of genus $1$ (corresponding to the term
$a_1$);
\item stable graphs with a single circuit, and all of whose vertices have
genus $0$ --- we call such a graph a \emph{necklace}.
\end{enumerate}
The presence of a logarithm in the term which contributes the necklaces is
related to the fact that there are $(n-1)!$ cyclic orders of $n$ objects.

The first few coefficients of $b_1$ are also easily calculated:
$$\begin{tabular}{|C|L|} \hline
n & \M v_{1,n} \\[2pt] \hline
1 & v_{1,1} + \half v_{0,3} \\[5pt]
2 & v_{1,2} + v_{1,1}v_{0,3} + \half\bigl( v_{0,4} + v_{0,3}^2 \bigr)
\\[5pt]
3 & v_{1,3} + 3v_{1,2}v_{0,3} + v_{1,1}v_{0,4} + \half\bigl( v_{0,5} +
3 v_{0,4}v_{0,3} + 2v_{0,3}^3 \bigr) \\[5pt]
4 & v_{1,4} + 6v_{1,3}v_{0,3} + 3v_{1,2}v_{0,4} + 15v_{1,2}v_{0,3}^2 +
v_{1,1}v_{0,5} \\
& \quad {}+ \half\bigl( v_{0,6} + 4v_{0,5}v_{0,3} + 3v_{0,4}^2 +
12v_{0,4}v_{0,3}^2 + 6v_{0,3}^4\bigr) \\[3pt] \hline
\end{tabular}$$

\section{The semi-classical approximation for modular operads}

In the theory of modular operads, one replaces the sequence of coefficients
$\{v_{g,n}\}$ considered above by a stable $\SS$-module, that is, a
sequence of $\SS_n$-modules $\V\(g,n\)$. The analogue of \eqref{M} is the
functor on stable $\SS$-modules which sends $\V$ to
\begin{equation} \label{MM}
\MM\V\(g,n\) = \colim_{G\in\Gamma_{g,n}} \bigotimes_{v\in\VERT(G)}
\V\(g(v),n(v)\) .
\end{equation}
Thus, the coefficients in \eqref{M} are promoted to vector spaces, the
product to a tensor product, the sum over stable graphs to a direct sum,
and the weight $|\Aut(G)|^{-1}$ to $\colim_{\Aut(G)}$, that is, the
coinvariants with respect to the finite group $\Aut(G)$. Note that this
definition makes sense in any symmetric monoidal category $\Cat$ with
finite colimits. We will need the Peter-Weyl theorem to hold for actions of
the symmetric group $\SS_n$ on $\Cat$; thus, we will suppose that $\Cat$ is
additive over a ring of characteristic zero.

\begin{definition}
The characteristic $\ch_n(\V)$ of an $\SS_n$-module is defined by the
formula
$$
\ch_n(\V) = \frac{1}{n!} \sum_{\sigma\in\SS_n} \Tr_\sigma(\V) p_\sigma \in
\Lambda_n\o K_0(\Cat) ,
$$
where $p_\sigma$ is the product of power sums $p_{|\CO|}$ over the orbits
$\CO$ of $\sigma$.
\end{definition}

Although this definition appears to require rational coefficients, this is
an artifact of the use of the power sums $p_n$; it is shown in \cite{I}
that the characteristic is a symmetric function of degree $n$ with values
in the Grothendieck group of the additive category $\Cat$. If
$\rk:\Lambda\to\Q[x]$ is the homomorphism defined by $h_n\mapsto x^n/n!$,
we have
$$
\rk(\ch_n(\V)) = [\V]/n!\in K_0(\Cat)\o\Q .
$$
Note that $\rk(f)$ is obtained from $f$ by setting the powers sums $p_n$ to
$0$ if $n>1$, and to $x$ if $n=1$.

The place of the generating functions $a_g$ and $b_g$ is now taken by
\begin{align*}
\AA_g & = \sum_{2(g-1)+n>0} \ch_n(\V\(g,n\)) \in \Lambda \ohat K_0(\Cat) , \\
\BB_g & = \sum_{2(g-1)+n>0} \ch_n(\MM\V\(g,n\)) \in \Lambda \ohat K_0(\Cat) .
\end{align*}
Theorem (8.13) of \cite{modular}, whose statement we now recall, calculates
$\BB_g$ in terms of $\AA_h$, $h\le g$. Let $\Delta$ be the ``Laplacian'' on
$\Lambda\(\hbar\)$ given by the formula
$$
\Delta = \sum_{n=1}^\infty \hbar^n
\left( \frac{n}{2} \frac{\p^2}{\p p_n^2} + \frac{\p}{\p p_{2n}} \right) .
$$
\begin{theorem} \label{modular}
If $\V$ is a stable $\SS$-module, then
$$
\sum_{g=0}^\infty \BB_g \hbar^{g-1} = \Log \biggl( \exp(\Delta) \Exp\Bigl(
\sum_{g=0}^\infty \AA_g \hbar^{g-1} \Bigr) \biggr) .
$$
\end{theorem}

There is also a formula for $\BB_0$ in terms of $\AA_0$. To state it, we
must recall the definition of the Legendre transform for symmetric
functions. Let
$$
\Lambda_*\ohat K_0(\Cat) = \{ f \in \Lambda\ohat K_0(\Cat) \mid \rk(f) =
x^2/2+O(x^3) \} .
$$
If $f$ is a symmetric function, let $f'=\p f/\p p_1$; this operation may be
expressed more invariantly as $p_1^\perp$ (Ex.\ I.5.3, Macdonald
\cite{Macdonald}).
\begin{definition}
The Legendre transform $\CL$ is the involution of $\Lambda_*\ohat
K_0(\Cat)$ characterized by the formula $(\CL f)\circ f' + f = p_1 f'$.
\end{definition}

The Legendre transform $\CL f$ of a function $f$ is characterized by the
formula $(\CL f)'\circ f'=x$. For symmetric functions, although the
analogue of this formula holds, in the form
$$
(\CL f)'\circ f' = h_1 ,
$$
the situation is not as simple, since there is no single notion of integral
for symmetric functions (the ``constant'' term may be any function of the
power sums $p_n$, $n>1$). Neverthless, there is a simple algorithm for
calculating $\CL f$ from $f$. Denote by $f_n$ and $g_n$ the coefficents of
$f$ and $g=\CL f$ lying in $\Lambda_n\o K_0(\Cat)$.
\begin{enumerate}
\item The formula $f'\circ(\CL f)'=h_1$ may be rewritten as
$$
\sum_{n=3}^N g_n' + \sum_{n=3}^N f_n' \circ
\Bigl( h_1 + \sum_{k=3}^{N-1} g_k' \Bigr) \cong 0 \mod{\Lambda_N\o K_0(\Cat)} .
$$
This gives a recursive procedure for calculating $g_n'$.
\item Having determined $g'$, we obtain $g$ from the formula $f=\CL g$, or
$g = p_1 g' - f\circ g'$.
\end{enumerate}

We now recall Theorem (7.17) of \cite{modular}, which is the generalization
to modular operads of Theorem \ref{legendre}.
\begin{theorem} \label{Legendre}
The symmetric function $h_2+\BB_0$ is the Legendre transform of
$e_2-\AA_0$.
\end{theorem}

The main result of this paper is a formula for $\BB_1$ in terms of
$\AA_1$ and $\AA_0$, generalizing Theorem \ref{semi}. If $f$ is a
symmetric function, write $\dot{f}=\p f/\p p_2=\half p_2^\perp f$.
\begin{theorem} \label{Semi}
$$
\BB_1 = \biggl( \AA_1 - \frac{1}{2} \sum_{n=1}^\infty
\frac{\phi(n)}{n} \log(1-\psi_n(\AA_0'')) +
\frac{\dot{\AA}_0(\dot{\AA}_0+1)}{1-\psi_2(\AA_0'')} \biggr) \circ
(h_1+\BB_0')
$$
Here, $\phi(n)$ is Euler's function, the number of prime residues modulo
$n$.
\end{theorem}

\begin{remark}
The first two terms inside the parentheses on the right-hand side of
Theorem \ref{Semi} are analogues of the corresponding terms in the formula
of Theorem \ref{semi}. In particular, the second of these terms is closely
related to the sum over necklaces in the definition of $\MM\V\(1,n\)$, as
is seem from the formula
$$
\sum_{n=1}^\infty \ch_n\bigl( \Ind_{\Z_n}^{\SS_n} \1 \bigr)
= - \sum_{n=1}^\infty \frac{\phi(n)}{n} \log(1-p_n) .
$$
The remaining term may be understood as a correction term, which takes into
account the fact that necklaces of $1$ or $2$ vertices have non-trivial
involutions (while those with more vertices do not). A proof of the theorem
could no doubt be given using this observation; however, we prefer to
derive it directly from Theorem \ref{modular}.

If we take the plethysm on the right of the formula of Theorem \ref{Semi}
with the symmetric function $h_1-\AA_0'$, and apply the formula
$(h_1+\BB_0')\circ(h_1-\AA_0')=h_1$, we obtain the equivalent formulation of
this theorem:
$$
\BB_1 \circ (h_1-\AA_0') = \AA_1 - \frac{1}{2} \sum_{n=1}^\infty
\frac{\phi(n)}{n} \log(1-\psi_n(\AA_0'')) +
\frac{\dot{\AA}_0(\dot{\AA}_0+1)}{1-\psi_2(\AA_0'')} .
$$
\end{remark}

\begin{proof}[Proof of Theorem \ref{Semi}]
The symmetric function $\BB_1$ is a sum over graphs obtained by attaching
forests whose vertices have genus $0$ to either a vertex of genus $1$, or
to a necklace. In other words,
$$
\BB_1 = \bigl( \AA_1 + \text{sum over necklaces} \bigr) \circ
(h_1+\BB_0') .
$$
To prove the theorem, we must calculate the sum over necklaces.

To do this, observe that a necklace is a graph with flags coloured red or
blue, such that each vertex has exactly two red flags, each edge is red,
and all tails are blue. Let $\w\(n\)$, $n\ge1$, be the sequence of
representations of $\SS_2\times\SS_n$
$$
\w\(n\) = \Res^{\SS_{n+2}}_{\SS_n\times\SS_2} \V\(0,n+2\) ;
$$
think of the first factor of the product $\SS_n\times\SS_2$ as acting on
the blue flags at a vertex, and the second factor as acting on the red
flags. Applying Theorem \ref{modular}, we see that
$$
\Log \bigl( \exp(1\o\Delta) \Exp(\Ch(\w)) \bigr) \in
\Lambda\ohat\Lambda\ohat K_0(\Cat)
$$
is the sum over stable graphs all of whose edges are red. To impose the
condition that all tails are blue, we set the variables $q_n$ to zero
before taking the Logarithm.

We now proceed to the explicit calculation. We set $\hbar=1$, since it
plays no r\^ole when all graphs have genus $1$. In writing elements of
$\Lambda\ohat\Lambda$, we will denote power sums in the first factor of
$\Lambda$ by $p_n$, and in the second by $q_n$.
\begin{lemma}
The characteristic $\Ch(\w)$ of $\w$ is the ``bisymmetric'' function
$$
\Ch(\w) = \half \AA_0''q_1^2 + \dot{\AA}_0 q_2 \in
\Lambda\ohat\Lambda_2\ohat K_0(\Cat) .
$$
\end{lemma}
\begin{proof}
We have $\Ch(\w)=h_2^\perp \AA_0\o h_2+e_2^\perp \AA_0\o e_2$. Expressing
this in terms of power sums, we have
\begin{align*}
h_2^\perp \AA_0\o h_2+e_2^\perp \AA_0\o e_2 &=
\bigl(\half(p_1^\perp)^2+p_2^\perp\bigr)\AA_0\o\half(q_1^2+q_2) +
\bigl(\half(p_1^\perp)^2-p_2^\perp\bigr)\AA_0\o\half(q_1^2-q_2) \\
&= \half (p_1^\perp)^2\AA_0\o q_1^2 + p_2^\perp \AA_0\o q_2 .
\qed\end{align*}
\def\qed{}
\end{proof}

From this lemma, it follows that
$$
\Exp\bigl( \Ch(\w) \bigr) = \prod_{n=1}^\infty \exp\Bigl( \psi_n(\AA_0'')
\frac{q_n^2}{2n} \Bigr) \prod_{n=1}^\infty \exp\Bigl( \psi_n(\dot{\AA}_0)
\frac{q_{2n}}{n} \Bigr) \in \Lambda\ohat\Lambda\ohat K_0(\Cat) ,
$$
We now apply the heat kernel and separate variables:
\begin{multline*}
\exp(1\o\Delta) \Exp\bigl( \Ch(\w) \bigr)\big|_{q_n=0}
= \prod_{\text{$n$ odd}} \exp\left( \frac{n}{2} \frac{\p^2}{\p q_n^2} \right)
\exp\Bigl( \psi_n(\AA_0'') \frac{q_n^2}{2n} \Bigr) \Big|_{q_n=0} \\
{} \times \prod_{\text{$n$ even}}
\exp\left( \frac{n}{2} \frac{\p^2}{\p q_n^2} + \frac{\p}{\p q_n} \right)
\exp\Bigl( \psi_n(\AA_0'') \frac{q_n^2}{2n} +
\psi_{n/2}(\dot{\AA}_0) \frac{2q_n}{n} \Bigr) \Big|_{q_n=0} .
\end{multline*}
We now insert the explicit formulas for the heat kernel of the Laplacian,
namely
$$
\exp\left( \frac{n}{2} \frac{\p^2}{\p q_n^2} \right) f(q_n) \big|_{q_n=0} =
\int_{-\infty}^\infty f(q_n) \exp\biggl( - \frac{q^2}{2n} \biggr)
\frac{dq}{\sqrt{2\pi n}} .
$$
For the odd variables, matters are quite straightforward:
\begin{align*}
\exp\left( \frac{n}{2} \frac{\p^2}{\p q_n^2} \right) \exp\Bigl(
\frac{q_n^2}{2n} \psi_n(\AA_0'') \Bigr) \Big|_{q_n=0} &=
\int_{-\infty}^\infty \exp\biggl( \psi_n(\AA_0'') \frac{q_n^2}{2n} -
\frac{q_n^2}{2n} \biggr) \frac{dq_n}{\sqrt{2\pi n}} \\
&= \bigl( 1 - \psi_n(\AA_0'') \bigr)^{-1/2} .
\end{align*}
For the even variables, things become a little more involved:
\begin{multline*}
\exp\left( \frac{n}{2} \frac{\p^2}{\p q_n^2} + \frac{\p}{\p q_n} \right)
\exp\Bigl( \psi_n(\AA_0'') \frac{q_n^2}{2n} +
\psi_{n/2}(\dot{\AA}_0) \frac{2q_n}{n} \Bigr) \Big|_{q_n=0} \\
\begin{aligned}
{} &= \exp\left( \frac{n}{2} \frac{\p^2}{\p q_n^2} \right) \exp\Bigl(
\psi_n(\AA_0'') \frac{q_n^2}{2n} + \psi_{n/2}(\dot{\AA}_0) \frac{2q_n}{n}
\Bigr) \Big|_{q_n=1} \\
{} &= \int_{-\infty}^\infty
\exp\biggl( \psi_n(\AA_0'') \frac{q_n^2}{2n} + \psi_{n/2}(\dot{\AA}_0)
\frac{2q_n}{n} - \frac{(q_n-1)^2}{2n} \biggr) \frac{dq_n}{\sqrt{2\pi n}} .
\end{aligned}
\end{multline*}
To perform this gaussian integral, we complete the square in the exponent:
\begin{multline*}
\psi_n(\AA_0'') \frac{q_n^2}{2n} + \psi_{n/2}(\dot{\AA}_0) \frac{2q_n}{n}
- \frac{(q_n-1)^2}{2n} \\
\begin{aligned}
{} &= - \bigl( 1-\psi_n(\AA_0'') \bigr) \frac{q_n^2}{2n} +
\bigl( 2\psi_{n/2}(\dot{\AA}_0) + 1 \bigr) \frac{q_n}{n} - \frac{1}{2n} \\
{} & = - \frac{1-\psi_n(\AA_0'')}{2n} \biggl( q_n -
\frac{2\psi_{n/2}(\dot{\AA}_0)+1}{1-\psi_n(\AA_0'')} \biggr)^2
+ \frac{2}{n} \Biggl( \frac{\psi_{n/2}(\dot{\AA}_0)
\bigl(\psi_{n/2}(\dot{\AA}_0)+1\bigr)}{1-\psi_n(\AA_0'')} \Biggr) .
\end{aligned}
\end{multline*}
Thus, the gaussian integral equals
$$
\bigl( 1-\psi_n(\AA_0'') \bigr)^{-1/2}
\exp\frac{2}{n} \Biggl(
\frac{\psi_{n/2}(\dot{\AA}_0)\bigl(\psi_{n/2}(\dot{\AA}_0)+1\bigr)}
{1-\psi_n(\AA_0'')} \Biggr) .
$$

Putting these calculations together, we see that
\begin{align*}
\exp(1\o\Delta) \Exp\bigl( \Ch(\w) \bigr) |_{q_n=0} &= \prod_{n=1}^\infty
\bigl( 1 - \psi_n(\AA_0'') \bigr)^{-1/2}
\exp\frac{1}{n} \Biggl(
\frac{\psi_n(\dot{\AA}_0)\bigl(\psi_n(\dot{\AA}_0)+1\bigr)}
{1-\psi_{2n}(\AA_0'')} \Biggr) \\
&= \prod_{n=1}^\infty \bigl( 1 - \psi_n(\AA_0'') \bigr)^{-1/2} \Exp \biggl(
\frac{\dot{\AA}_0 (\dot{\AA}_0+1)}{1-\psi_2(\AA_0'')} \biggr) ,
\end{align*}
and, applying the operation $\Log$, that
$$
\Log \bigl( \exp(1\o\Delta) \Exp\bigl( \Ch(\w) \bigr) |_{q_n=0} \bigr) =
\Log \prod_{n=1}^\infty \bigl( 1 - \psi_n(\AA_0'') \bigr)^{-1/2} +
\frac{\dot{\AA}_0(\dot{\AA}_0+1)}{1-\psi_2(\AA_0'')} .
$$
The proof of the theorem is completed by the following lemma., applied to
$f=1-\AA_0''$.
\begin{lemma}
Let $f\in\Lambda\ohat K_0(\Cat)$ have constant term equal to $1$; that is,
$\rk(f)=1+O(x)$. Then
$$
\Log \prod_{n=1}^\infty \psi_n(f)^{-1/2} = - \frac{1}{2} \sum_{n=1}^\infty
\frac{\phi(n)}{n} \log(\psi_n(f)) .
$$
\end{lemma}
\begin{proof}
By definition,
$$
\Log \prod_{n=1}^\infty \psi_n(f)^{-1/2} = \sum_{k=1}^\infty
\frac{\mu(k)}{k} \log \prod_{n=1}^\infty \psi_{nk}(f)^{-1/2} = -
\frac{1}{2} \prod_{n=1}^\infty \Bigl( \sum_{d|n} \frac{\mu(d)}{d} \Bigr)
\log(\psi_n(f)) .
$$
The lemma follows from the formula
$$
\sum_{d|n} \frac{\mu(d)}{d} = \frac{\phi(n)}{n} ,
$$
which follows by M\"obius inversion from $\sum_{d|n}\phi(d)=n$.
\end{proof}
\def\qed{}
\end{proof}

\begin{corollary} \label{sEMI}
Define $a_g=\rk(\AA_g)$, $b_g=\rk(\BB_g)$, and $\dot{a}_0=\rk(\dot{\AA}_0)$. Then
we have
$$
a_1 \circ (x-a_0') = a_1 - \half \log(1-a_0'') + \dot{a}_0
(\dot{a}_0+1) .
$$
\end{corollary}

\begin{example}
Suppose $\V\(0,n\)=\1$ is the trivial one-dimensional representation for
all $n\ge3$, while $\V\(1,n\)=0$. Then $\MM\V\(1,n\)$ is an $\SS_n$-module
whose rank is the number of graphs in $\Gamma^0_{1,n}$, where
$\Gamma^0_{1,n}\subset\Gamma_{1,n}$ is the subset of stable graphs all of
whose vertices have genus $0$. We have
$$
\AA_0 = \sum_{n=3}^\infty h_n = \exp\Bigl( \sum_{n=1}^\infty \frac{p_n}{n}
\Bigr) - 1 - h_1 - h_2 .
$$
Theorem \ref{Semi} leads to the following results; the calculations were
performed using J.~Stembridge's symmetric function package \texttt{SF} for
\texttt{maple} \cite{SF}.
$$\begin{tabular}{|C|L|L|} \hline
n & \ch_n\bigl(\MM\V\(1,n\)\bigr) & |\Gamma_{1,n}^0| \\ \hline
1 & s_{1} & 1 \\[5pt]
2 & 3\,s_{2} & 3 \\[5pt]
3 & 7\,s_{3}+4\,s_{21} & 15 \\[5pt]
4 & 20\,s_{4}+17\,s_{31}+14\,s_{2^2}+4\,s_{21^2} & 111 \\[5pt]
5 & 52\,s_{5}+78\,s_{41}+71\,s_{32}+33\,s_{31^2}+34\,s_{2^21}+4\,s_{21^3}+s_{1^5} &
1104 \\ \hline
\end{tabular}$$

An explicit formula for the generating function of the numbers
$|\Gamma_{1,n}^0|$ may be obtained from Corollary \ref{sEMI}, using the
formulas $a_0'=e^x-1-x$, $a_0''=e^x-1$ and $\dot{a}_0=\half(e^x-1)$.
\begin{proposition}
$$
\sum_{n=1}^\infty |\Gamma^0_{1,n}| \frac{x^n}{n!} = \Bigl( - \frac{1}{2}
\log \bigl( 2 - e^x \bigr) + \frac{1}{4} (e^{2x}-1) \Bigr) \circ
(1+2x-e^x)^{-1} .
$$
\end{proposition}
\end{example}

\section{The $\SS_n$-equivariant Hodge polynomial of $\Mbar_{1,n}$}

A more interesting application of Theorem \ref{Semi} is to the stable
$\SS$-module in the category of $\Z$-graded mixed Hodge structures
$$
\V\(g,n\) = H^\bull_c(\CM_{g,n},\C) .
$$

Let $\KHM$ be the Grothendieck group of mixed Hodge structures. The
$\SS_n$-equivariant Serre polynomial $\Serre^{\SS_n}(\CM_{g,n})$ is by
definition the characteristic $\ch_n(\V\(g,n\))\in\Lambda_n\o\KHM$. It
follows from the usual properties of Serre polynomials (see \cite{I} or
Proposition (6.11) of \cite{modular}) that $\ch_n(\MM\V\(g,n\))$ is the
$\SS_n$-equivariant Serre polynomial $\Serre^{\SS_n}(\Mbar_{g,n})$ of the
moduli space $\Mbar_{g,n}$ of stable curves. Since the moduli space
$\Mbar_{g,n}$ is a complete smooth Deligne-Mumford stack, its $k$th
cohomology group carries a pure Hodge structure of weight $k$; thus, the
Hodge polynomial of $\Mbar_{g,n}$ may be extracted from
$\Serre^{\SS_n}(\Mbar_{g,n})$. Using Theorem \ref{Semi}, we will calculate
the Serre polynomials $\Serre^{\SS_n}(\Mbar_{1,n})$.

It is shown in \cite{gravity} (see also \cite{I}) that
$$
\AA_0 = \sum_{n=3}^\infty \Serre^{\SS_n}(\CM_{0,n}) = \frac{\displaystyle
\biggl\{ \prod_{n=1}^\infty
(1+p_n)^{\frac{1}{n}\sum_{d|n}\mu(n/d)(1+\LL^d)} \biggr\} - 1}{\LL^3-\LL} -
\frac{h_1}{\LL^2-\LL} - \frac{h_2}{\LL+1} ,
$$
where $\LL$ is the pure Hodge structure $\C(-1)$ of weight $2$. Theorem
\ref{Legendre} implies that
$$
h_2 + \BB_0 = h_2 + \sum_{n=3}^\infty \Serre^{\SS_n}(\Mbar_{0,n})
$$
is the Legendre transform of $e_2-\AA_0$; this was used in \cite{gravity} to
calculate $\Serre^{\SS_n}(\Mbar_{0,n})$.

Let $\SSS_{2k+2}$ be the pure Hodge structure
$\operatorname{gr}^W_{2k+1}H^1_c(\CM_{1,1},\operatorname{Sym}^{2k}\HH)$,
where $\HH$ is the local system $R^1\pi_*\Q$ of rank $2$ over the moduli
stack of elliptic curves. (Here, $\pi:\Mbar_{1,2}\to\Mbar_{1,1}$ is the
universal elliptic curve.) This Hodge structure has the following
properties:
\begin{enumerate}
\item $\SSS_{2k+2}=F^0\SSS_{2k+2}\oplus\overline{F^0\SSS_{2k+2}}$;
\item there is a natural isomorphism between $F^0\SSS_{2k+2}$ and the space
of cusp forms $S_{2k+2}$ for the full modular group $\SL(2,\Z)$. (In
particular, $\SSS_{2k+2}=0$ for $k\le4$.)
\end{enumerate}

It is shown in \cite{II} that
\begin{multline*}
\AA_1 = \sum_{n=1}^\infty \Serre^{\SS_n}(\CM_{1,n}) = \res_0 \Biggl[ \left(
\frac{\prod_{n=1}^\infty
(1+p_n)^{\frac{1}{n}\sum_{d|n}\mu(n/d)(1-\om^d-\LL^d/\om^d+\LL^d)} - 1}
{1-\om-\LL/\om+\LL} \right) \\ \times \left( \sum_{k=1}^\infty \biggl(
\frac{\SSS_{2k+2}+1}{\LL^{2k+1}} \biggr) \om^{2k} - 1 \right) \bigl(
\om-\LL/\om \bigr) d\om \Biggr] ,
\end{multline*}
where $\res_0[\alpha]$ is the residue of the one-form $\alpha$ at the
origin.

We may now apply Theorem \ref{Semi} to calculate the generating function of
the $\SS_n$-equivariant Serre polynomials $\Serre^{\SS_n}(\Mbar_{1,n})$. We
do not give the details, since they are quite straightforward, though the
resulting formulas are tremendously complicated when written out in
full. However, we do present some sample calculations, performed with the
package \texttt{SF}.
$$\begin{tabular}{|C|L|L|} \hline
n & \Serre\bigl(\Mbar_{1,n}\bigr) & \chi(\Mbar_{1,n}) \\[2pt] \hline
1 & (\LL+1)s_1 & 2 \\[5pt]
2 & (\LL^2+2\LL+1)s_2 & 4 \\[5pt]
3 & (\LL^3+3\LL^2+3\LL+1)s_3+(\LL^2+\LL)s_{21} & 12 \\[5pt]
4 & (\LL^4+4\LL^3+7\LL^2+4\LL+1)s_4+(2\LL^3+4\LL^2+2\LL)s_{31}
+(\LL^3+2\LL^2+\LL)s_{2^2} & 49 \\[5pt]
5 & (\LL^5+5\LL^4+12\LL^3+12\LL^2+5\LL+1)s_5
+(3\LL^4+11\LL^3+11\LL^2+3\LL)s_{41} & 260 \\
& {}+(2\LL^4+7\LL^3+7\LL^2+2\LL)s_{32}+(\LL^3+\LL^2)(s_{31^2}+s_{2^21}) &
\\[2pt] \hline
\end{tabular}$$

In a table at the end of the paper, we give a table of non-equivariant
Serre polynomials of $\Mbar_{1,n}$ for $n\le15$; these give an idea of the
way in which the Hodge structures $\SSS_{2k+2}$ typically enter into the
cohomology. In particular, we see that the even-dimensional cohomology of
the moduli spaces $\Mbar_{1,n}$ is spanned by Hodge structures of the form
$\Q(\ell)$, while the odd dimensional cohomology is spanned by Hodge
structures of the form $\SSS_{2k+2}(\ell)$.

The rational cohomology groups of $\Mbar_{1,n}$ satisfy Poincar\'e duality:
there is a non-degenerate $\SS_n$-equivariant pairing of Hodge structures
$$
H^k(\Mbar_{1,n},\Q) \o H^{2n-k}(\Mbar_{1,n},\Q) \to \Q(-n) .
$$
Unfortunately, our formula for $\Serre^{\SS_n}(\Mbar_{1,n})$ does not
render this duality manifest.

\section{The Euler characteristic of $\Mbar_{1,n}$}

As an application of Corollary \ref{sEMI}, we give an explicit formula for
the generating function of the Euler characteristics $\chi(\Mbar_{1,n})$.
\begin{theorem} \label{funny}
Let $g(x)\in x+x^2\Q\[x\]$ be the solution of the equation
$$
2g(x)-(1+g(x))\log(1+g(x))=x .
$$
Then
$$
\sum_{n=1}^\infty \chi(\Mbar_{1,n}) \frac{x^n}{n!} = - \frac{1}{12}
\log\bigl(1+g(x)\bigr) - \frac{1}{2} \log\bigl(1-\log(1+g(x))\bigr) +
\epsilon(g(x)) ,
$$
where
$$
\epsilon(x) = \frac{1}{12} \bigl( 19\,x + 23\,x^2/2 + 10\,x^3/3 + x^4/2
\bigr) .
$$
\end{theorem}
\begin{proof}
We apply Corollary \ref{sEMI} with the data
\begin{align*}
a_0' &= \sum_{n=2}^\infty \chi(\CM_{0,n+1}) \frac{x^n}{n!} =
\sum_{n=2}^\infty (-1)^n (n-2)! \frac{x^n}{n!} = (1+x)\log(1+x)-x , \\
a_0'' &= \log(1+x) , \quad \dot{a}_0 = \frac{1}{4} x(x+2) , \\
a_1 &= \chi(\CM_{1,1}) x + \chi(\CM_{1,2}) \frac{x^2}{2} +
\chi(\CM_{1,3}) \frac{x^3}{6} + \chi(\CM_{1,4}) \frac{x^4}{24}
+ \frac{1}{12} \sum_{n=5}^\infty (-1)^n (n-1)! \frac{x^n}{n!} \\
&= x + \frac{x^2}{2} - \frac{1}{12} \log(1+x) + \frac{1}{12} \Bigl( x +
\frac{x^2}{2} + \frac{x^3}{3} - \frac{x^4}{4} \Bigr) , \\
\end{align*}
where we have used that $\chi(\CM_{1,1})=\chi(\CM_{1,2})=1$ and
$\chi(\CM_{1,3})=\chi(\CM_{1,4})=0$. The function $g(x)$ of the statement
of the theorem is $x+\BB_0'(x)$.
\end{proof}

The following corollary was shown us by D. Zagier.
\begin{corollary} \label{don}
$$
\chi(\Mbar_{1,n}) \sim \frac{(n-1)!}{4(e-2)^n} \Bigl( 1 + C
n^{-1/2} + O\bigl(n^{-3/2}\bigr) \Bigr) ,
$$
where
$$
C = \sqrt{\frac{e-2}{18\pi e}} ( 1 + 4e + 9e^2 + 4e^3 + 2e^4 ) \approx
18.31398807 .
$$
\end{corollary}
\begin{proof}
To show this, we analytically continue $g(x)$ to the domain
$\C\setminus[e-2,\infty)$. The resulting function has an asymptotic
expansion of the form
$$
g(x) \sim e - 1 - \sqrt{2e(e-2-x)} + \sum_{k=3}^\infty a_k (e-2-x)^{k/2} .
$$
The asymptotics \eqref{don} follow by applying Cauchy's integral formula to
the right-hand side of Theorem \ref{funny}, with contour the circle
$|x|=e-2$.
\end{proof}

The peculiar polynomial $\epsilon(x)$ of Theorem \ref{funny} combines the
error terms in the formula for $\chi(\CM_{1,n})$ with the correction terms
involving $\dot{a}_0$ in Corollary \ref{sEMI}. Omitting the term
$\epsilon(g(x))$ in Theorem \ref{funny}, we obtain the generating function
not of the Euler characteristics $\chi(\Mbar_{1,n})$, but rather of the
virtual Euler characteristics $\chi_v(\Mbar_{1,n})$ of the underlying
smooth moduli stack (orbifold). The asymptotic behaviour of the virtual
Euler characteristics is the same as that of the Euler characteristics,
with $C$ replaced by $\widetilde{C}=\bigl(\frac{e-2}{18\pi e}\bigr)^{1/2}
\approx 0.06835794$. The ratio between these Euler characteristics has the
asymptotic behaviour
$$
\frac{\chi(\Mbar_{1,n})}{\chi_v(\Mbar_{1,n})} \sim (C-\widetilde{C})
n^{-1/2} + O(n^{-1}) ,
$$
giving a statistical measure of the ramification of $\Mbar_{1,n}$ for large
$n$.

\begin{sideways}
$$\begin{tabular}{|C|L|} \hline
n & \Serre(\Mbar_{1,n}) \\ \hline
1 & \LL+1 \\
2 & \LL^2+2\,\LL+1 \\
3 & \LL^3+5\,\LL^2+5\,\LL+1 \\
4 & \LL^4+12\,\LL^3+23\,\LL^2+12\,\LL+1 \\
5 & \LL^5+27\,\LL^4+102\,\LL^3+102\,\LL^2+27\,\LL+1 \\
6 & \LL^6+58\,\LL^5+421\,\LL^4+756\,\LL^3+421\,\LL^2+58\,\LL+1 \\
7 & \LL^7+121\,\LL^6+1612\,\LL^5+5077\,\LL^4+5077\,\LL^3+1612\,\LL^2+12\,\LL+1 \\
8 & \LL^8+248\,\LL^7+5802\,\LL^6+31072\,\LL^5+52402\,\LL^4+31072\,\LL^3
+5802\,\LL^2+248\,\LL+1 \\
9 & \LL^9+503\,\LL^8+19925\,\LL^7+175036\,\LL^6+480097\,\LL^5+480097\,\LL^4
+175036\,\LL^3+19925\,\LL^2+503\,\LL+1 \\
10 & \LL^{10}+1014\,\LL^9+66090\,\LL^8+920263\,\LL^7+3975949\,\LL^6
+6349238\,\LL^5+3975949\,\LL^4+920263\,\LL^3+66090\,\LL^2+1014\,\LL+1 \\
11 & \LL^{11}+2037\,\LL^{10}+213677\,\LL^9+4577630\,\LL^8+30215924\,\LL^7
+74269967\,\LL^6+30215924\,\LL^5+\ldots+1 - \SSS_{12} \\
12 & \LL^{12}+4084\,\LL^{11}+677881\,\LL^{10}+21793602\,\LL^9
+213725387\,\LL^8+784457251\,\LL^7+1196288936\,\LL^6+\ldots+4084\,\LL+1
-11(\LL+1)\SSS_{12} \\
13 & \LL^{13}+8179\,\LL^{12}+2120432\,\LL^{11}+100226258\,\LL^{10}
+1424858788\,\LL^9+7603002045\,\LL^8+17095248952\,\LL^7+\ldots \\
& \quad {} - (66\,\LL^2+429\,\LL+66)\SSS_{12} \\
14 & \LL^{14}+16370\,\LL^{13}+6563147\,\LL^{12}+448463866\,\LL^{11}
+ 9049174765\,\LL^{10}+68547770726\,\LL^9
+221071720149\,\LL^8+324314241400\,\LL^7+\ldots \\
& \quad {} - (286\,\LL^3+6006\,\LL^2+286\,\LL)\SSS_{12} \\
15 & \LL^{15}+32753\,\LL^{14}+20153930\,\LL^{13}+1963368663\,\LL^{12}
+55228789080\,\LL^{11} + 581636563570\,\LL^{10}
+2627427327522\,\LL^9+5488190927216\,\LL^8+\ldots \\
& \quad {} -
(1001\,\LL^4+53053\,\LL^3+186263\,\LL^2+53053\LL+1001)\SSS_{12}-\SSS_{16}
\\ \hline
\end{tabular}$$
\end{sideways}

\end{document}